\documentclass[10pt,draft,twoside,leqno,a4paper]{article}
\usepackage{amsmath,amsthm,amsfonts,amssymb,verbatim}     
\newtheorem{Theorem}{Theorem}[section]    
\newtheorem{Lemma}{Lemma}[section]    
\newtheorem{Corollary}{Corollary}[section]    
\newtheorem{Proposition}{Proposition}[section]    
\newtheorem{Remark}{Remark}[section]   
   
\newcommand{\R}{\mathbb R}

\newcommand{\us}{{\widetilde u}}   
\newcommand{\uph}{U} 
\newcommand{\vph}{V}   
\newcommand{\sph}{S}
\newcommand{\qph}{Q}   
\newcommand{\oly}{\overline y} 
\newcommand{\uly}{\underline y} 
\newcommand{\xk}{{X^k}}   
\begin{document}
\title{On a nonlinear elliptic system from 
Maxwell-Chern-Simons vortex theory}        
\author{Tonia Ricciardi\thanks{
Partially supported by PRIN 2000 ``Variational Methods and Nonlinear 
Differential Equations"}\\
{\small Dipartimento di Matematica e Applicazioni}\\
{\small Universit\`a di Napoli Federico II}\\
{\small Via Cintia}\\
{\small 80126 Naples, Italy}\\
{\small fax: +39 081 675665}\\
{\small e-mail: tonia.ricciardi@unina.it}
}  
\date{July 30, 2002}      
\maketitle
\begin{abstract}
We define an abstract nonlinear elliptic system,
admitting a variational structure, and
including the vortex equations
for some Maxwell-Chern-Simons gauge theories
as special cases. We analyze the asymptotic
behavior of its solutions, and we provide a general simplified framework for the
asymptotics previously derived in those special cases.
As a byproduct of our abstract formulation, 
we also find some new qualitative properties of solutions.
\end{abstract}
\begin{tabbing}
\={\sc key words}: \=nonlinear elliptic system,
Chern-Simons vortex theory\\
\>{\sc MCS 2000 subject classification}: 35J60
\end{tabbing}
\setcounter{section}{-1} 
\section{Introduction} 
\label{sec:introduction} 
Motivated by the analysis of vortex configurations for
the self-dual $U(1)$ Max\-well-Chern-Simons model introduced in \cite{LLM}
(see also Yang~\cite{Y},  Dun\-ne~\cite{D}, Jaffe and Taubes~\cite{JT}), 
we considered in \cite{R} solutions $(u,v)$ for the system:
\begin{align}
\label{mcsu1a'}
-\Delta u=&q(v-e^u)-4\pi\sum_{j=1}^n\delta_{p_j}
&&\text{on}\ \Sigma\\
\label{mcsu1b'}
-\Delta v=&q\left\{e^u(1-v)-q(v-e^u)\right\}
&&\text{on}\ \Sigma,
\end{align} 
where $\Sigma$ is a compact Riemannian 2-manifold without boundary,
$n\ge0$ is an integer, $p_j\in\Sigma$ for $j=1,\ldots,n$,
$\Delta$ denotes the Laplace-Beltrami operator
and $q>0$ is a constant. 
It is of both mathematical and physical interest to understand the asymptotic behavior
of solutions to \eqref{mcsu1a'}--\eqref{mcsu1b'}
as $q\to+\infty$.
In \cite{R} we provided a rigorous proof of the formal
asymptotics derived in \cite{LLM}, in any relevant norm.
Namely, we showed that if $(u,v)$ are (distributional)
solutions for \eqref{mcsu1a'}--\eqref{mcsu1b'}
with $q\to+\infty$, then there exists a solution $u_\infty$ for the equation
\begin{equation}
\label{csu1}
-\Delta u_\infty=e^{u_\infty}(1-e^{u_\infty})-4\pi\sum_{j=1}^n\delta_{p_j}
\qquad\qquad\text{on}\ \Sigma
\end{equation}
such that $(e^u,v)\to(e^{u_\infty},e^{u_\infty})$ in 
$C^h(\Sigma)\times C^h(\Sigma)$, for any $h\ge0$.
(Note that $e^u$, $e^{u_\infty}$ are smooth).
Such a result completed our previous convergence result obtained
with Tarantello~\cite{RT}, 
where the asymptotics for $v$ was established in
the $L^2$-sense only. See also Chae and Kim~\cite{CK}. 

More recently, Chae and Nam~\cite{CN} analyzed
an elliptic system, whose solutions describe vortex configurations 
for the self-dual $CP(1)$ Maxwell-Chern-Simons model introduced 
in \cite{KLL}. Their system (in a special case) is given by:
\begin{align}
\label{mcscp1a'}
\Delta\uph=&2\qph(-\vph+\sph-\frac{1-e^\uph}{1+e^\uph})+4\pi\sum_{j=1}^n\delta_{p_j}
&&\text{on}\ \Sigma\\
\label{mcscp1b'}
\Delta\vph=&-4\qph^2(-\vph+\sph-\frac{1-e^\uph}{1+e^\uph})+\qph\frac{4e^\uph}{(1+e^\uph)^2}\vph
&&\text{on}\ \Sigma
\end{align}
where $\Sigma$ and $p_1,\ldots,p_n$ are as in \eqref{mcsu1a'}--\eqref{mcsu1b'},
$U,V$ are the unknown functions and $S\in\R$, $Q>0$ are given constants.
Among other results, they derive an asymptotic behavior 
as $Q\to+\infty$ analogous to that of system
\eqref{mcsu1a'}--\eqref{mcsu1b'}.\par
Thus it is natural to seek a common underlying structure
for \eqref{mcsu1a'}--\eqref{mcsu1b'} and \eqref{mcscp1a'}--\eqref{mcscp1b'},
which would allow such asymptotic behaviors.\par
Our aim in this note is to identify a general nonlinear system 
including 
\eqref{mcsu1a'}--\eqref{mcsu1b'} and \eqref{mcscp1a'}--\eqref{mcscp1b'}
as special cases, and to
show that the asymptotic behaviors
described above are in fact a {\em general property} of its 
solutions. We believe that our proof of the asymptotics 
for our abstract system simplifies 
the previous approaches. 
Our formulation will also
allow us to find some new qualitative properties of solutions.\par
More precisely, we consider (distributional) solutions $(\us,v)$
for the system:
\begin{align}
\label{genmcsa'}
&-\Delta\us=q(v-f(e^\us))-4\pi\sum_{j=1}^n\delta_{p_j}
&&\text{on}\ \Sigma\\
\label{genmcsb'}
&-\Delta v=q\left[f'(e^\us)e^\us(s-v)-q(v-f(e^\us))\right]
&&\text{on}\ \Sigma.
\end{align}
Here $\Sigma$ and $p_1,\ldots,p_n$ are as in \eqref{mcsu1a'}--\eqref{mcsu1b'},
$f=f(t)$, $t\ge0$ is smooth and strictly increasing, $s\in\R$
satisfies
$f(0)<s<\sup_{t>0}f(t)$.
We shall later show that when $n=0$,
system \eqref{genmcsa'}--\eqref{genmcsb'}
only admits the trivial solution $(f(e^\us),v)=(s,s)$.
Without loss of generality, we assume $\text{vol}\Sigma=1$.\par
Clearly, when $f(t)=t$ and $s=1$, system \eqref{genmcsa'}--\eqref{genmcsb'}
reduces to \eqref{mcsu1a'}--\eqref{mcsu1b'}.
On the other hand, setting $v:=\vph-\sph$, $s:=-\sph$, $q:=2\qph$,
system \eqref{genmcsa'}--\eqref{genmcsb'} reduces to system \eqref{mcscp1a'}--\eqref{mcscp1b'}
with $f$ defined by $f(t)=(t-1)/(t+1)$.\par
As already mentioned, we are interested in the asymptotic behavior of
solutions as $q\to+\infty$.
By a formal analysis of \eqref{genmcsa'}--\eqref{genmcsb'} 
we expect that, up to subsequences, $(\us,v)$ converges
to some solution $(\widetilde u_\infty,f(e^{\widetilde u_\infty}))$,
for the equation:
\begin{equation}
\label{gencs}
-\Delta\widetilde u_\infty=f'(e^{\widetilde u_\infty})e^{\widetilde u_\infty}
(s-f(e^{\widetilde u_\infty}))-4\pi\sum_{j=1}^n\delta_{p_j}
\qquad\qquad\text{on}\ \Sigma.
\end{equation}
Our main result states that this is indeed the case,
with respect to any relevant norm:
\begin{Theorem}
\label{thm:main}
Let $(\us,v)$ be (distributional) solutions to \eqref{genmcsa'}--\eqref{genmcsb'},
with $q\to+\infty$. There exists a (distributional) solution $\widetilde u_\infty$ to
\eqref{gencs}
such that a subsequence, still denoted $(\us,v)$, satisfies:
\begin{align*}
(e^\us,v)&\to\left(e^{\widetilde u_\infty},f(e^{\widetilde u_\infty})\right)\qquad\qquad
\text{in}\ C^h(\Sigma)\times C^h(\Sigma),\ \forall h\ge0.\\
\end{align*}
\end{Theorem}
In order to work in suitable Sobolev spaces, 
it is standard (see \cite{Y}) to define
a ``Green's function" $u_0$, solution for the problem
\begin{align*}
&-\Delta u_0=4\pi\left(n-\sum_{j=1}^n\delta_{p_j}\right)
\qquad\text{on}\ \Sigma\\
&\int_\Sigma u_0=0
\end{align*}
(see \cite{Aubin} for the unique existence of $u_0$).
Setting $\us=u_0+u$, we obtain the equivalent system
for $(u,v)\in H^1(\Sigma)\times H^1(\Sigma)$:
\begin{align}
\label{genmcsa}
&-\Delta u=q\left(v-f(e^{u_0+u})\right)-4\pi n
&&\text{on}\ \Sigma\\
\label{genmcsb}
&-\Delta v=
q\left[f'(e^{u_0+u})e^{u_0+u}(s-v)-q\left(v-f(e^{u_0+u})\right)\right]
&&\text{on}\ \Sigma,
\end{align}
where $e^{u_0}$ is {\em smooth}.
We also note that system \eqref{genmcsa}--\eqref{genmcsb}
admits a {\em variational formulation}. Indeed, 
solutions $(u,v)$ to \eqref{genmcsa}--\eqref{genmcsb}
correspond to critical points $v\in H^2(\Sigma)$
for the functional:
\begin{align*}
I(u)=\frac{1}{2q^2}&\int(\Delta u)^2+\frac{1}{2}\int|\nabla u|^2\\
+&\frac{1}{q}\int f'(e^{u_0+u})e^{u_0+u}|\nabla(u_0+u)|^2
+\frac{1}{2}\int(f(e^{u_0+u})-s)^2+4\pi n\int u.
\end{align*}
Since
\begin{equation*}
e^{u_0}|\nabla u_0|^2=\Delta e^{u_0}+4\pi n,
\end{equation*}
the function $e^{u_0}|\nabla u_0|^2$ is smooth, and 
$I$ is well-defined on $H^2(\Sigma)$ by Sobolev embeddings.
To see how critical points for $I$ correspond to
solutions for \eqref{genmcsa}--\eqref{genmcsb},
we solve \eqref{genmcsa} for $v$:
\[
v=q^{-1}(-\Delta u+4\pi n)+f(e^{u_0+u}).
\]
Substituting into \eqref{genmcsb}, we obtain
the fourth-order equation:
\begin{align}
\nonumber
\frac{1}{q^2}\Delta^2 u-\Delta u
-&\frac{1}{q}\left[\Delta f(e^{u_0+u})+f'(e^{u_0+u})e^{u_0+u}\Delta(u_0+u)\right]\\
\label{fourthorder}
+&f'(e^{u_0+u})e^{u_0+u}(f(e^{u_0+u})-s)+4\pi n=0
\qquad\text{on}\ \Sigma.
\end{align}
Integration by parts shows that
\begin{align*}
\frac{d}{dt}\big|_{t=0}
\int f'(e^{u_0+u+t\phi})&e^{u_0+u+t\phi}|\nabla(u_0+u+t\phi)|^2\\
=&-\int\left\{\Delta f(e^{u_0+u})+f'(e^{u_0+u})e^{u_0+u}\Delta(u_0+u)\right\},
\end{align*}
and thus critical points for $I$ correspond to solutions
for \eqref{fourthorder}.
We shall exploit this variational structure in order to
study the {\em multiplicity} of solutions to \eqref{genmcsa'}--\eqref{genmcsb'}
in a forthcoming note.\par
The remaining part of this note is devoted to the proof of 
Theorem~\ref{thm:main}.
The main point of the proof is to obtain
{\em a priori estimates} for $\us-u_0$ and for $v$ independent of $q\to+\infty$ in
the Sobolev spaces $H^k$, for every $k\ge0$.
To this end, in Section~\ref{sec:apriori} we first establish 
some preliminary estimates
in $L^\infty$ and $H^1$.
In Section~\ref{sec:iteration}, exploiting the {\em specific structure}
of system \eqref{genmcsa}--\eqref{genmcsb}, we set up
an iteration in the framework of the Banach algebras
$H^k\cap L^\infty$, for $k\ge0$, which yields the desired estimates.\par
Henceforth we denote by $C>0$ a general constant
independent of $q$, which may vary from line to line.
Unless otherwise specified, 
all equations are defined on $\Sigma$ and
all integrals are taken over $\Sigma$
with respect to the Lebesgue measure.
\section{A priori estimates}
\label{sec:apriori}
Our aim in this section is to 
establish estimates in $H^1$ and in $L^\infty$ 
for $e^\us$ and $v$, as stated in the following
\begin{Proposition}
\label{prop:basis}
There exists a constant $C>0$ independent of $q\to+\infty$,
such that:
\begin{equation}
\|e^\us\|_{H^1\cap L^\infty}+\|v\|_{H^1\cap L^\infty}
+\|q(v-f(e^\us))\|_{L^2}\le C.
\end{equation}
\end{Proposition}
We shall first obtain some pointwise estimates,
which depend 
on the increasing monotonicity of $f$ in an essential way:
\begin{Lemma}
\label{lem:pointwise}
The following estimates hold, pointwise on $\Sigma$:
\begin{align*}
\tag{i}&f(0)\le f(e^{\us})\le s\\
\tag{ii}&f(0)\le v\le s.
\end{align*}
\end{Lemma}
\begin{Corollary}
\label{cor:novortex}
If $n=0$, then $(e^\us,v)=(f^{-1}(s),s)$.
\begin{proof}
Suppose $n=0$. 
Integrating \eqref{genmcsa'} and \eqref{genmcsb'} we find that
\[
q\int(v-f(e^u))=0=\int f'(e^u)e^u(s-v).
\]
By Lemma \ref{lem:pointwise}--(ii) we have $v\le s$.
Since $f'(e^u)e^u>0$, the above identity implies $v\equiv s$.
Then $\Delta v\equiv0$ and thus \eqref{genmcsb'} implies $q(s-f(e^u))\equiv0$,
that is, $f(e^u)\equiv s$, as asserted.
\end{proof}
\end{Corollary}
As a consequence of Lemma \ref{lem:pointwise}, the nonlinearity $f$
may be {\em truncated}. Therefore in what follows, without loss of generality,
we assume that:
\begin{equation}
\label{ftruncation}
\sup_{t>0}\{|f(t)|+|f'(t)|+|f''(t)|\}\le C.
\end{equation}
\begin{proof}[Proof of Lemma \ref{lem:pointwise}]
Let $\bar x\in\Sigma$ be such that $\us(\bar x)=\max_\Sigma\us$.
Then $\overline x\neq p_j$ for all $j=1,\ldots,n$
and \eqref{genmcsa'} implies that
\begin{equation*}
f(e^{\us(\bar x)})\le v(\bar x).
\end{equation*}
Now we equivalently rewrite equation \eqref{genmcsb'} in the form:
\begin{align}
\label{vequiv}
-\Delta v+q^2\left(1+\frac{1}{q}f'(e^\us)e^\us\right)v
=q^2\left(f(e^\us)+\frac{s}{q}f'(e^\us)e^\us\right).
\end{align}
Let $\oly,\uly\in\Sigma$ such that
$v(\oly)=\max_\Sigma v$, $v(\uly)=\min_\Sigma v$.
Then, the maximum principle applied to \eqref{vequiv}
implies that
\begin{align}
\label{maxprinctov}
\frac{f(e^{\us(\uly)})+\frac{s}{q}f'(e^{\us(\uly)})e^{\us(\uly)}}
{1+\frac{1}{q}f'(e^{\us(\uly)})e^{\us(\uly)}}
\le v
\le\frac{f(e^{\us(\oly)})+\frac{s}{q}f'(e^{\us(\oly)})e^{\us(\oly)}}
{1+\frac{1}{q}f'(e^{\us(\oly)})e^{\us(\oly)}},
\end{align}
pointwise on $\Sigma$.
If $\oly=p_j$ for some $j=1,\ldots,n$, then $e^{\us(\oly)}=0$
and therefore the second inequality in \eqref{maxprinctov} implies:
$v(\oly)\le f(0)$.
Since $f(0)<s$ by assumption, the second part of (ii) is established in 
this case.
(In fact, we can show that $\oly\neq p_j$, for all $j=1,\ldots,n$,
see Remark \ref{rem:qualitative} below). 
If $\oly\neq p_j$ for all $j=1,\ldots,n$, then
we observe that by increasing monotonicity of $f$ we have:
\[
f(e^{\us(\oly)})
\le f(e^{\us(\bar x)})\le v(\bar x)\le v(\oly).
\]
Inserting into the second inequality in \eqref{maxprinctov}, we derive:
\[
\left(1+\frac{1}{q}f'(e^{\us(\oly)})e^{\us(\oly)}\right)v(\oly)
\le v(\oly)+\frac{s}{q}f'(e^{\us(\oly)})e^{\us(\oly)},
\]
which, recalling that $f'>0$, in turn yields :
\begin{equation}
e^{\us(\oly)}v(\oly)\le s\,e^{\us(\oly)},
\end{equation}
with $e^{\us(\oly)}>0$. 
Hence (i) and the second part of (ii) follow.
It remains to show that $v\ge f(0)$.
By (i), we know that $s-f(e^{\us(\uly)})\ge0$.
Therefore, by the increasing monotonicity of $f$:
\begin{equation}
\label{forvlowerbound}
\frac{s-f(e^{\us(\uly)})}{1+\frac{1}{q}f'(e^{\us(\uly)})}
\le s-f(e^{\us(\uly)})\le s-f(0).
\end{equation}
Consequently, combining the first inequality in \eqref{maxprinctov}
and \eqref{forvlowerbound}, we obtain:
\begin{align*}
v(\uly)\ge
\frac{f(e^{\us(\uly)})+\frac{s}{q}f'(e^{\us(\uly)})e^{\us(\uly)}}
{1+\frac{1}{q}f'(e^{\us(\uly)})e^{\us(\uly)}}
=s-\frac{s-f(e^{\us(\uly)})}{1+\frac{1}{q}f'(e^{\us(\uly)})}
\ge f(0),
\end{align*}
and the proof of (ii) is complete.
\end{proof}
As already mentioned in the proof of Lemma~\ref{lem:pointwise}, we can actually show
that $v$ does not attain its maximum at $p_j$, $j=1,\ldots,n$:
\begin{Remark}
\label{rem:qualitative}
If $v$ is constant, then $n=0$ and $(e^\us,v)=(f(s),s)$.
In particular, if $n>0$, then $v$ cannot be a constant.
Furthermore, if $n>0$, then $v$ attains its maximum on $\Sigma\setminus\{p_1,\ldots,p_n\}$.
\end{Remark}
\begin{proof}
Suppose $v\equiv k=\text{constant}$.
Then integrating \eqref{genmcsa'} and \eqref{genmcsb'} we obtain
\[
\int f'(e^\us)e^\us(s-k)=q\int(k-f(e^\us))=4\pi n
\]
and therefore
\[
k=s-\frac{4\pi n}{\int f'(e^\us)e^\us}\le s.
\]
If $k=s$, then $n=0$ and by Corollary~\ref{cor:novortex} we have $k=s=f^{-1}(e^\us)$. 
Thus, the statement of the 
lemma is established in this case.\par
Therefore we assume $k<s$. In particular $n>0$, and thus $\us$
is not constant. Furthermore, $t=e^\us$ attains values in 
$[0,\delta)$ for some $\delta>0$.
Setting $\phi(\us)=f(e^\us)$, we have from \eqref{genmcsb'} that
$\phi$ satisfies the differential equation
\[
(s-k)\phi'=q(k-\phi)\qquad\qquad,\us\le-M
\]
for some $M>0$ and thus
\[
\phi(\us)=Ce^{-q\us/(s-k)}+k.
\]
Recalling the definition of $\phi$, it follows that $f$
has the form
\[
f(t)=Ct^{-q/(s-k)}+k,
\]
which is singular at $t=0$, contradiction.
Now, if $v$ attains its maximum at some $p_j$, then by \eqref{maxprinctov} necessarily
$v\equiv f(0)$. Hence, $n=0$.
\end{proof}
The next estimate will be used to derive $H^1$-bounds for $e^\us$
and $v$:
\begin{Lemma}
\label{lem:graduestimate}
We have:
\begin{equation*}
\int e^\us|\nabla\us|^2\le C.
\end{equation*}
\end{Lemma}
\begin{proof}
Multiplying equation \eqref{genmcsa'} by $e^\us$
and integrating by parts,
we obtain 
\[
q\int e^{\us}(v-f(e^\us))
=\int e^{\us}|\nabla\us|^2\ge0.
\]
By the pointwise estimates in Lemma \ref{lem:pointwise},
it follows that:
\begin{equation}
\label{graduprelim}
\frac{1}{q}\int e^{\us}|\nabla\us|^2\le C.
\end{equation}
Multiplying \eqref{genmcsb'} by $e^\us$ and integrating, we find
\begin{align}
\label{intbyparts}
q\int e^\us(v-f(e^\us))=\int e^{2\us}f'(e^\us)(s-v)
+\frac{1}{q}\int e^\us\Delta v.
\end{align}
Integration by parts yields:
\[
\frac{1}{q}\int e^\us\Delta v
=-\int ve^\us(v-f(e^\us))
+\frac{1}{q}\int ve^\us|\nabla\us|^2.
\]
Hence, by the pointwise estimates as in Lemma \ref{lem:pointwise},
and taking into account \eqref{graduprelim}, we conclude that
\[
\frac{1}{q}\left|\int e^\us\Delta v\right|\le C.
\]
Inserting into \eqref{intbyparts}, recalling Lemma \ref{lem:pointwise},
we derive that
\[
q\int e^\us(v-f(e^\us))\le C,
\]
and thus it follows that
\begin{equation*}
\int e^{\us}|\nabla\us|^2=q\int e^{\us}(v-f(e^\us))
\le C.
\end{equation*}
\end{proof}
The next identity is the main step in deriving
the $H^1$-estimate for $v$ and the $L^2$-estimate for
$q(v-f(e^\us))$:
\begin{Lemma}
\label{lem:identity}
The following identity holds:
\begin{align}
\label{mainidentity}
\int|\nabla v|^2+q^2\int(v-f(e^\us))^2
=\int(s-v)\left(f''(e^\us)e^\us+f'(e^\us)\right)e^\us|\nabla\us|^2.
\end{align}
\end{Lemma}
\begin{proof}
We compute:
\[
\Delta f(e^\us)
=\left(f''(e^\us)e^\us+f'(e^\us)\right)e^\us|\nabla\us|^2
+f'(e^\us)e^\us\Delta\us.
\]
Therefore $f(e^\us)$ satisfies the equation:
\begin{align}
\label{feuequation}
-\Delta f(e^u)+qf'(e^\us)e^\us\,f(e^\us)
=qf'(e^\us)e^\us v-\left(f''(e^\us)e^\us
+f'(e^\us)\right)e^\us|\nabla\us|^2.
\end{align}
Integrating \eqref{feuequation}, we obtain 
\begin{equation}
\label{integralidentity1}
q\int f'(e^\us)e^\us(v-f(e^\us))
=\int\left(f''(e^\us)e^\us+f'(e^\us)\right)e^\us|\nabla\us|^2
\end{equation}
Now we multiply \eqref{genmcsb'} by $v-f(e^\us)$ and integrate to obtain:
\begin{align*}
\int-\Delta v(v-f(e^\us))
=q\int f'(e^\us)e^\us(s-v)(v-f(e^\us))-q^2\int(v-f(e^\us))^2.
\end{align*}
Integrating by parts and using \eqref{feuequation}
we find:
\begin{align*}
\int-\Delta&v(v-f(e^\us))
=\int|\nabla v|^2+\int v\Delta f(e^\us)\\
=&\int|\nabla v|^2-q\int vf'(e^\us)e^\us(v-f(e^\us))
+\int v\left(f''(e^\us)e^\us+f'(e^\us)\right)e^\us|\nabla\us|^2.
\end{align*}
Equating left hand sides in the last two identities, we obtain 
\begin{align*}
&\int|\nabla v|^2+q^2\int(v-f(e^\us))^2
+\int v\left(f''(e^\us)e^\us+f'(e^\us)\right)e^\us|\nabla\us|^2\\
&\qquad\qquad\qquad\qquad=sq\int f'(e^\us)e^\us(v-f(e^\us)),
\end{align*}
and thus identity \eqref{mainidentity} is established.
\end{proof}
Now we can finally provide the
\begin{proof}[Proof of Proposition \ref{prop:basis}]
Lemma~\ref{lem:pointwise} readily implies
$\|e^\us\|_{L^\infty}\le C$ and $\|v\|_{L^\infty}\le C$.
In order to obtain the $H^1$-estimate for $e^\us$,
it suffices to observe that by Lemma~\ref{lem:pointwise}--(i)
and by Lemma~\ref{lem:graduestimate} we have:
\[
\int|\nabla e^\us|^2=\int e^{2\us}|\nabla\us|^2
\le C\int e^\us|\nabla\us|^2\le C.
\]
Therefore, we are left to estimate $\|\nabla v\|_{L^2}$
and $\|q(v-f(e^\us))\|_{L^2}$.
Using identity \eqref{mainidentity}, we have:
\begin{align*}
\int|\nabla v|^2+q^2\int(v-f(e^\us))^2
\le&\|s-v\|_\infty\|f''(e^\us)e^\us+f^(\us)\|_\infty\int e^\us|\nabla\us|^2\\
\le&C\int e^\us|\nabla\us|^2\le C,
\end{align*}
where we again used Lemma~\ref{lem:pointwise}
and Lemma~\ref{lem:graduestimate} in order to derive the last step.
\end{proof}
\section{Iteration}
\label{sec:iteration}
The aim of this section is to obtain bounds for solutions
in $H^k$, for every $k\ge0$, as given in the following
\begin{Proposition}
\label{prop:iteration}
For all $k\ge0$ there exists a constant $C>0$
(possibly depending on $k$)
such that:
\begin{equation*}
\|\us-u_0\|_{H^k}+\|v\|_{H^k}\le C.
\end{equation*}
\end{Proposition}
It will be convenient to define the Banach spaces
$X^0:=L^\infty$, $\xk:=H^k\cap L^\infty$ for $k\ge1$,
endowed with the norms $\| \cdot\|_\xk:=\|\cdot\|_{H^k}+\|\cdot\|_{L^\infty}$.
We recall the well-known Sobolev-Gagliardo-Nirenberg inequality,
see e.g.\ \cite{Nirenberg}:
\begin{equation*}
\|D^j u\|_{L^p(\R^n)}\le C\|D^k u\|_{L^r(\R^n)}^a\|u\|_{L^q(\R^n)}^{1-a}
\qquad\qquad\forall u\in C_c^\infty(\R^n),
\end{equation*}
where
\begin{align*}
&\frac{1}{p}=\frac{j}{k}+a\left(\frac{1}{r}-\frac{m}{n}\right)
+(1-a)\frac{1}{q}\\
&\frac{j}{k}\le a\le 1.
\end{align*}
Taking $k=2$, $a=j/k$, $q=\infty$, and using partitions of unity
on $\Sigma$, we obtain:
\begin{equation}
\label{xkproperty}
\|D^ju\|_{L^{2k/j}}\le C\|D^k u\|_{L^2}^{j/k}\|u\|_{L^\infty}^{1-j/k}
\qquad\qquad\forall u\in C^\infty(\Sigma).
\end{equation}
By \eqref{xkproperty} and the H\"older inequality
that if $u_1,\ldots,u_t\in\xk$ and $\beta_1,\ldots,\beta_t$
are multi-indices such that $|\beta_1|+\cdots+|\beta_t|=k$,
then the product $D^{\beta_1}u_1\cdots D^{\beta_t}u_t\in L^2$
and
\[
\|D^{\beta_1}u_1\cdots D^{\beta_t}u_t\|_{L^2}\le C\|u_1\|_\xk\cdots\|u_t\|_\xk.
\] 
In particular, $\xk$ is a Banach algebra for every $k\ge0$, 
i.e.,
\begin{equation*}
\|u_1u_2\|_\xk\le C\|u_1\|_\xk\|u_2\|_\xk.
\end{equation*}
We shall need the following
\begin{Lemma}
\label{lem:nonlinearitybound}
Let $F\in C^\infty(\Sigma\times\R)$,
$G\in C^\infty(\Sigma\times\R\times\R\times\R^2)$.
Then for all $k\ge0$ there exists constants $C_k=C_k(\|u\|_{L^\infty})$,
$C_k'=C_k'(\|u\|_{L^\infty},\|v\|_{L^\infty},\|\nabla u\|_{L^\infty})$, such that:
\begin{align*}
&\|F(x,u)\|_\xk\le C_k(1+\|u\|_\xk^k)\\
&\|G(x,u,v,\nabla u)\|_{X^{k-1}}\le C_k'(1+\|u\|_\xk^{k-1}+\|v\|_{X^{k-1}}).
\end{align*} 
\end{Lemma}
\begin{proof}
Denote by $\alpha$ a multi-index such that $|\alpha|=k$.
It suffices to observe that
\[
D^\alpha F(x,u)=\sum_{|\alpha_1|+\cdots+|\alpha_h|=|\alpha|}
F^{(h)}(x,u)D^{\alpha_1}u\cdots D^{\alpha_h}u
\]
and therefore, recalling \eqref{xkproperty}, we have
\begin{align*}
\|D^\alpha F(x,u)\|_{L^2}\le&C(\|u\|_{L^\infty})
\sum_{|\alpha_1|+\cdots+|\alpha_h|=|\alpha|}
\|D^{\alpha_1}u\cdots D^{\alpha_h}u\|_{L^2}\\
\le&C(\|u\|_{L^\infty})\sum_{|\alpha_1|+\cdots+|\alpha_h|=|\alpha|}
\|D^{\alpha_1}u\|_{L^{2|\alpha|/|\alpha_1|}}
\cdots\|D^{\alpha_h}u\|_{L^{2|\alpha|/|\alpha_h|}}\\
\le&C(\|u\|_{L^\infty})(1+\|u\|_\xk^k).
\end{align*}
Since obviously $\|F(x,u)\|_{L^\infty}\le C(\|u\|_{L^\infty})$,
the first estimate is established.
The estimate for $G(x,u,v,\nabla u)$ is obtained analogously.
\end{proof}
Now we observe that \eqref{genmcsb} is of the form:
\begin{equation}
\label{abstracteq}
-\Delta u+q^2(1+\frac{1}{q}c)u=q^2f.
\end{equation}
We shall need some a priori estimates for solutions to \eqref{abstracteq}.
The next two results state that, under suitable assumptions,
a solution $u$  for \eqref{abstracteq}
satisfies the same regularity properties as the right hand side $f$,
{\em independently} of $q\to+\infty$.
\begin{Lemma}
\label{lem:lpestimate}
Suppose $u$ is a solution for \eqref{abstracteq} with $c\in L^\infty$
and $f\in L^p$ for some $1\le p\le+\infty$.
Then there exist $q_k>0$ and $C>0$ independent of $u$ such that
\begin{equation}
\label{lpestimate}
\|u\|_{L^p}\le C\|f\|_{L^p},
\end{equation}
for all $q\ge q_k$.
\end{Lemma}
\begin{proof}
For $p=+\infty$, the statement follows by the maximum principle:
\begin{equation}
\label{abstractmax}
\|u\|_{L^\infty}\le\left\|\frac{f}{1+\frac{1}{q}c}\right\|_{L^\infty}
\le\frac{\|f\|_{L^\infty}}{1-\frac{1}{q}\|c\|_{L^\infty}},
\end{equation}
hence for large $q$ we find:
\[
\|u\|_{L^\infty}\le C\|f\|_{L^\infty}.
\]
Now we assume $2\le p<+\infty$.
Multiplying \eqref{abstracteq} by $|u|^{p-2}u$
and integrating by parts, we find
\begin{align*}
(p-1)\int|u|^{p-2}|\nabla u|^2+\int(1+\frac{1}{q}c)|u|^p
=\int f|u|^{p-2}u
\end{align*}
It follows that
\begin{align*}
\|u\|_{L^p}^p\le\frac{\int|f||u|^{p-1}}{1-\frac{1}{q}\|c\|_{L^\infty}}
\end{align*}
and therefore, by the H\"older inequality,
\[
\|u\|_{L^p}\le\frac{\|f\|_{L^p}}{1-\frac{1}{q}\|c\|_{L^\infty}},
\]
and hence the asserted estimate is established in the case $2\le p\le+\infty$.
In the remaining case $1\le p<2$, we proceed by duality.
Let $\varphi$ be defined by
\[
-\Delta\varphi+(1+\frac{1}{q}c)\varphi=|u|^{p-2}u.
\]
Then \eqref{lpestimate} with $p\ge2$ yields $\|\varphi\|_{L^{p'}}\le C\|u\|_{L^p}^{p-1}$.
Multiplying \eqref{abstracteq} by $\varphi$ and integrating, we find:
\begin{align*}
\int|u|^p
=\int-\Delta u\varphi+\int(1+\frac{1}{q}c)u\varphi=\int f\varphi.
\end{align*}
Consequently,
\[
\int|u|^p\le\|f\|_{L^p}\|\varphi\|_{L^{p'}}\le C\|f\|_{L^p}\|u\|_{L^p}^{p-1},
\]
and the asserted estimate is established also for $1\le p<2$.
\end{proof}
\begin{Lemma}
\label{lem:ellipticproperty}
Let $c,f\in\xk$ and suppose that $u$ satisfies:
\eqref{abstracteq}.
For every $k\ge0$ there exist $q_k>0$, $C_k>0$
such that
\begin{equation*}
\|u\|_\xk\le C_k\|f\|_\xk,
\end{equation*}
for all $q\ge q_k$.
\end{Lemma}
\begin{proof}
Denote by $\alpha$ a multi-index, $|\alpha|=k$.
Multiplying \eqref{abstracteq} by $D^{2\alpha}u$
and integrating by parts, we obtain:
\begin{align*}
\int|\nabla D^\alpha u|^2+q^2\int D^\alpha[(1+\frac{1}{q}c)u]D^\alpha u=
q^2\int D^\alpha u D^\alpha f,
\end{align*}
Therefore, 
\begin{align*}
\int D^\alpha[(1+\frac{1}{q}c)u]D^\alpha u
\le\int D^\alpha f D^\alpha u,
\end{align*}
and thus we estimate:
\begin{align*}
\int(D^\alpha u)^2\le&\int D^\alpha f D^\alpha u
-\frac{1}{q}\int D^\alpha(cu)D^\alpha u\\
\le&\|D^\alpha f\|_{L^2}\|D^\alpha u\|_{L^2}
+\frac{1}{q}\|D^\alpha(cu)\|_{L^2}\|D^\alpha u\|_{L^2}\\
\le&\|u\|_\xk\|f\|_\xk+\frac{1}{q}\|cu\|_\xk\|u\|_\xk\\
\le&\|u\|_\xk\|f\|_\xk+\frac{1}{q}\|c\|_\xk\|u\|_\xk^2,
\end{align*}
where we have used \eqref{xkproperty} to derive the last line.
Since $\alpha$ is an arbitrary multi-index satisfying $|\alpha|=k$, 
we conclude from the above and \eqref{abstractmax} that
\[
\|u\|_\xk\le C(\|f\|_\xk+\frac{1}{q}\|c\|_\xk\|u\|_\xk).
\]
Now the asserted estimate follows easily.
\end{proof}
At this point, it is useful to note that 
$q(v-f(e^{u_0+u}))$
also satisfies an equation of the form \eqref{abstracteq}.
In fact, it is convenient to set
\[
w:=q(v-f(e^{u_0+u})) 
\]
and to consider $w$ as a third unknown function.
Then the triple $(u,v,w)$ satisfies a system of the following simple form:
\begin{align}
\label{iteru}
&-\Delta u=w-4\pi n\\
\label{iterv}
&-\Delta v+q^2[1+\frac{1}{q}c(x,u)]v
=q^2F_q(x,u)\\
\label{iterw}
&-\Delta w+q^2[1+\frac{1}{q}c(x,u)]w
=q^2G_q(x,u,v,\nabla u)
\end{align}
where
\begin{align*}
c(x,u)=&f'(e^{u_0+u})e^{u_0+u}\\
F_q(x,u)=&f(e^{u_0+u})+\frac{s}{q}f'(e^{u_0+u})e^{u_0+u}\\
G_q(x,u,v,\nabla u)
=&f'(e^{u_0+u})e^{u_0+u}(s-v)\\
&\qquad+\frac{1}{q}\left(f''(e^{u_0+u})e^{u_0+u}
+f'(e^{u_0+u})\right)e^{u_0+u}|\nabla({u_0+u})|^2.
\end{align*}
Proposition~\ref{prop:iteration} will follow by a bootstrap argument
applied to \eqref{iteru}--\eqref{iterv}--\eqref{iterw}.
In order to start the procedure, we need:
\begin{Lemma}
\label{lem:uest}
The following estimates hold:
\begin{align*}
\tag{i}
&\|u\|_{X^1}+\|v\|_{X^1}+\|w\|_{X^0}\le C\\
\tag{ii}
&\|u\|_{L^\infty}\le C\\
\tag{iii}
&\|\nabla u\|_{L^\infty}\le C
\end{align*}
\end{Lemma}
\begin{proof}
Proof of (i).
Multiplying \eqref{genmcsa} by $u-\int u$ and integrating,
we have:
\begin{align*}
\int|\nabla u|^2
=&q\int(v-f(e^\us))(u-\int u)\\
\le&\|q(v-f(e^\us))\|_2\|u-\int u\|_2\le C\|\nabla u\|_2,
\end{align*}
where the last inequality follows by Lemma \ref{lem:pointwise} and 
by the Poincar\'e inequality.
Hence 
$\|\nabla u\|_2\le C$.
By Lemma \ref{lem:pointwise}--(ii), 
we have that $e^\us\le C$, and thus we only have to show that
$\int u\ge-C$.
To this end, we first observe that integrating \eqref{genmcsa}
and \eqref{genmcsb} we obtain:
\[
\int f'(e^{u_0+u})e^{u_0+u}(s-v)=q\int(v-f(e^{u_0+u}))=4\pi n.
\]
On the other hand, we have in a straightforward manner:
\[
\int f'(e^{u_0+u})e^{u_0+u}(s-v)
\le C\int e^{u_0+u}
\le Ce^{\int u}\|e^{u_0}\|_\infty\int e^{u-\int u}\le C\int e^{u-\int u}.
\]
Hence, recalling the Moser-Trudinger inequality (see \cite{Aubin}) and the
estimate
for $\|\nabla u\|_2$, we conclude that
\begin{align*}
4\pi n\le Ce^{\int u}\int e^{u-\int u}
\le Ce^{\int u}e^{\gamma\int|\nabla u|^2}\le Ce^{\int u},
\end{align*}
which establishes (i).
Proof of (ii).
Since $\|w\|_{L^2}\le C$, by (i) and elliptic regularity
we obtain $\|u\|_{H^2}\le C$. Then Sobolev embeddings 
yield $\|\nabla u\|_{L^p}\le C$, for any $1\le p<+\infty$
and $\|u\|_{L^\infty}\le C$, which establishes (ii).
Proof of (iii).
By \eqref{iterw}, $\|\nabla u\|_{L^p}\le C$ and 
Lemma~\ref{lem:lpestimate} imply
that $\|w\|_{L^p}\le C$, for any $1\le p<+\infty$.
Then \eqref{iteru} and Sobolev embeddings 
yield $\|u\|_{W^{2,p}}\le C$, for any $1\le p<+\infty$.
For $p>2$, the Sobolev embeddings yield (iii).
\end{proof}
Now we can provide the
\begin{proof}[Proof of Proposition~\ref{prop:iteration}]
We argue by induction on $k\in\mathbb N_0$.\par
{\bf CLAIM A:} There holds:
\begin{equation*}
\|u\|_{X^1}+\|v\|_{X^1}+\|w\|_{X^0}\le C.
\end{equation*}
The above follows by Proposition~\ref{prop:basis}
and by Lemma~\ref{lem:uest}.\par
{\bf CLAIM B:} Suppose:
\begin{equation*}
\|u\|_\xk+\|v\|_\xk+\|w\|_{X^{k-1}}\le C_k.
\end{equation*}
Then:
\begin{equation*}
\|u\|_{X^{k+1}}+\|v\|_{X^{k+1}}+\|w\|_\xk\le C_{k+1}.
\end{equation*}
Indeed, 
\begin{align*}
\|w\|_{X^{k-1}}\le C
&\Rightarrow\|u\|_{X^{k+1}}\le C&&\text{by \eqref{iteru} and standard elliptic regularity}\\
&\Rightarrow\|v\|_{X^{k+1}}\le C
&&\text{by \eqref{iterv}, Lemma~\ref{lem:nonlinearitybound} and Lemma~\ref{lem:ellipticproperty}}\\
&\Rightarrow\|w\|_\xk\le C
&&\text{by \eqref{iterw}, Lemma~\ref{lem:nonlinearitybound} 
and Lemma~\ref{lem:ellipticproperty}}.
\end{align*}
Now Claim A, Claim B and a standard induction argument conclude the proof.
\end{proof}
Finally, we can prove our main result:
\begin{proof}[Proof of Theorem~\ref{thm:main}]
Let $(u,v)$ be solutions to system \eqref{genmcsa}--\eqref{genmcsb},
with $q\to+\infty$.
By the a priori estimates as stated in Proposition~\ref{prop:iteration}
and by standard compactness arguments,  
there exist $u_\infty$, $v_\infty$ such that up to subsequences
$u\to u_\infty$ and $v\to v_\infty$ in $C^h$, for all $h\ge0$.
We write \eqref{genmcsa} in the form: 
\[
v=f(e^{u_0+u})+\frac{1}{q}(-\Delta u+4\pi n).
\]
Taking limits, we find $v_\infty=f(e^{u_0+u_\infty})$.
Furthermore, taking limits in \eqref{genmcsb}, we obtain
\[
q(v-f(e^{u_0+u}))\to f'(e^{u_0+u_\infty})e^{u_0+u_\infty}(s-f(e^{u_0+u_\infty})),
\]
where the convergence holds in $C^h$, for any $h\ge0$.
Consequently, taking limits in \eqref{genmcsa}, we find that
$u_\infty$ satisfies:
\begin{equation}
-\Delta u_\infty=f'(e^{u_0+u_\infty})e^{u_0+u_\infty}(s-f(e^{u_0+u_\infty}))
-4\pi\sum_{j=1}^n\delta_{p_j}.
\end{equation}
Setting $\widetilde u_\infty=u_0+u_\infty$, we conclude the proof
of Theorem~\ref{thm:main}.
\end{proof}
\section*{Acknowledgements}
I am grateful to Professors Danielle Hilhorst, 
Masayasu Mimura and Gabriella Tarantello for interesting and stimulating discussions.

\vfill
\hfill http://cds.unina.it/$\widetilde\ $tonricci
\end{document}